\documentclass[12pt]{article}

\usepackage[margin=1in]{geometry}  
\usepackage{graphicx}              
\usepackage{amsmath}               
\usepackage{amsfonts}              
\usepackage{amsthm}                
\usepackage{mathabx}
\usepackage{hyperref}
\newtheorem{Proposition}{Proposition}
\newtheorem*{Remark}{Remark}

\theoremstyle{definition}

\begin{document}

\title{Quantum Lichnerowicz – Poisson complex}

\author{Valerii Sopin\\
  \small email: \texttt{vVs@myself.com}}
\maketitle

\begin{abstract}
Using the curved bc--beta-gamma system (a tensor product of a Heisenberg and a Clifford vertex algebra)  we introduce quantum analogy of Lichnerowicz differential. As follows we suggest new  machinery for finding the Lichnerowicz–Poisson cohomology groups for any Poisson manifold. Moreover, the defined provides new invariant.

\end{abstract}

\section{Introduction}

A Poisson manifold is a smooth n-dimensional manifold  endowed with a Poisson bivector field, viz., a skew-symmetric contravariant tensor $P$ of rank 2 satisfying the Jacobi identity
$$\sum_{l=1}^{n} (P^{lj}\frac{\partial P^{ik}}{\partial x_l}+P^{li}\frac{\partial P^{kj}}{\partial x_l} +P^{lk}\frac{\partial P^{ji}}{\partial x_l})=0,$$
where $P$ is locally given by $\sum\limits_{i<j} P^{ij}(x) \partial_{x_i} \wedge \partial_{x_j}$.

A Poisson structure $P$ on a manifold defines geometric object, the Lichnerowicz differential $d_{L}$ discovered in $[1]$. It acts on multivector fields by the formula $d_{L}:=[[P, -]],$ where $[[-, -]]$ denotes the canonical Schouten bracket.

We are going to give a quantum analogy of $d_{L}$ based on the  paper of Malikov, Schechtman and Vaintrob $[2]$, who introduced a sheaf of vertex superalgebras $\Omega^{ch}$ attached to any smooth variety M, called the chiral de Rham complex, which is used in understanding the “stringy” invariants, such as the elliptic genera. If M is n–dimensional, the fibers of $\Omega^{ch}$ are isomorphic as vertex superalgebras to a completion of the $bc-\beta\gamma$ system on $n$ generators, or in physics terminology, to the tensor product of the bosonic and fermionic ghost systems.  

We would like to outline here that the chiral de Rham complex is an example of the general localization pattern $[3][4]$. 

The paper $[5]$ served as a main motive for this research. We would like to shed the light on mathematical part of $[5]$. On the contemporary state-of-the art, the interested reader is referred to $[3-8]$ and references therein.

The corresponding cohomology of Lichnerowicz differential $d_{L}$ is called the Lichnerowicz – Poisson cohomology (or LP-cohomology). It is a useful tool in Poisson Geometry, as it provides framework to express deformation and quantization obstructions. For every smooth Poisson manifold there is a natural homomorphism from its de Rham cohomology to its Lichnerowicz–Poisson cohomology. For symplectic manifolds, this homomorphism is an isomorphism $[7]$. But, generally, the LP-cohomology space are very large and their structure is known only in some particular cases. The quantum version of $d_{L}$ helps to clarify this issue by encoding everything in OPEs (the operator product expansion).

\section{Quantum Lichnerowicz differential}

Consider  Poisson n-dimensional manifold M with a Poisson tensor $P = \sum\limits_{i<j}  P^{ij}(x) \partial_{x_i} \wedge \partial_{x_j},$ where $P^{ij}$ are analytic functions.

Let $\alpha^{i}=\gamma^i+b^i dz$, $\theta_{i}=c_i+\beta_i dz$, where $b_i, c_i, \gamma_i, \beta_i$ is corresponding $bc-\beta\gamma$ system on $n$ generators on $M$ according to $[2].$

Assign
$$\oint P^{ij}(\alpha) \theta_i\theta_j = \oint P^{ij}(\gamma+b dz) (c_i+\beta_i dz)(c_j+\beta_j dz).$$

The last is equal to $$\oint (P^{ij}(\gamma) (c_i\beta_j - \beta_i c_j) + \partial_k P^{ij}(\gamma)c_i c_j b^k) dz.$$

Define $$d_{qL}= \{\oint P^{ij}(\alpha) \theta_i\theta_j, -\}.$$

\begin{Proposition} $d^2_{qL}=0.$\end{Proposition}

\textit{Proof.}  The associated non-zero  Operator Product Expansions (OPEs) for $bc$ and $\beta\gamma$ systems:
$$b(z)c(w) \sim \frac{\hbar}{z-w}, c(z)b(w) \sim \frac{\hbar}{z-w};\; \beta(z)\gamma(w) \sim \frac{-\hbar}{z-w}, \gamma(z)\beta(w) \sim \frac{\hbar}{z-w}.$$

Further, without loss of generality $\hbar:=1.$

As $P^{ij}(\gamma)$ is analytic function, it can locally be written via a convergent power series. Moreover, notice $$(\gamma^k)^n(z) \beta_{\hat{k}}(w) \sim \delta_{k, \hat{k}}\frac{n (\gamma^k)^{n-1}}{z-w}.$$

Using Wick's theorem we obtain the following:

Term $P^{ij}(\gamma(z)) (c_i(z)\beta_j(z) - \beta_i(z) c_j(z)) P^{\hat{i}\hat{j}}(\gamma(w)) (c_{\hat{i}}(w)\beta_{\hat{j}}(w)- \beta_{\hat{i}}(w) c_{\hat{j}}(w))$ yields 
$$ \frac{1}{z-w} \partial_{\hat{j}}  P^{ij}(\gamma) P^{\hat{i}\hat{j}}(\gamma) c_i c_{\hat{i}} \beta_j  -  \frac{1}{z-w} \partial_{j} P^{\hat{i}\hat{j}}(\gamma) P^{ij}(\gamma) c_i c_{\hat{i}}\beta_{\hat{j}}  -\frac{1}{(z-w)^2}  \partial_{\hat{j}}  P^{ij}(\gamma) \partial_{j} P^{\hat{i}\hat{j}}(\gamma) c_i c_{\hat{i}}  -$$
$$- \frac{1}{z-w} \partial_{\hat{j}}  P^{ij}(\gamma) P^{\hat{i}\hat{j}}(\gamma) c_j c_{\hat{i}} \beta_i  + \frac{1}{z-w} \partial_{i} P^{\hat{i}\hat{j}}(\gamma) P^{ij}(\gamma) c_j c_{\hat{i}}\beta_{\hat{j}} +\frac{1}{(z-w)^2} \partial_{\hat{j}}  P^{ij}(\gamma) \partial_{i} P^{\hat{i}\hat{j}}(\gamma) c_j c_{\hat{i}} -$$

$$ - \frac{1}{z-w} \partial_{\hat{i}} P^{ij}(\gamma) P^{\hat{i}\hat{j}}(\gamma) c_i c_{\hat{j}} \beta_j  + \frac{1}{z-w} \partial_{j}P^{\hat{i}\hat{j}}(\gamma) P^{ij}(\gamma) c_i c_{\hat{j}}\beta_{\hat{i}} + \frac{1}{(z-w)^2} \partial_{\hat{i}}  P^{ij}(\gamma) \partial_{j}P^{\hat{i}\hat{j}}(\gamma) c_i c_{\hat{j}} +$$
$$ + \frac{1}{z-w} \partial_{\hat{i}}  P^{ij}(\gamma) P^{\hat{i}\hat{j}}(\gamma) c_j c_{\hat{j}} \beta_i  -  \frac{1}{z-w} \partial_{i} P^{\hat{i}\hat{j}}(\gamma) P^{ij}(\gamma) c_j c_{\hat{j}}\beta_{\hat{i}} - \frac{1}{(z-w)^2}  \partial_{\hat{i}}  P^{ij}(\gamma) \partial_{i} P^{\hat{i}\hat{j}}(\gamma) c_j c_{\hat{j}}.$$

Term $$P^{ij}(\gamma(z)) (c_i(z)\beta_j(z) - \beta_i(z) c_j(z)) \partial_{\hat{k}} P^{\hat{i}\hat{j}}(\gamma(w)) c_{\hat{i}}(w)c_{\hat{j}}(w) b^{\hat{k}}(w) +$$ $$+ \partial_k P^{ij}(\gamma(z))c_i(z) c_j(z) b^k(z) P^{\hat{i}\hat{j}}(\gamma(w)) (c_{\hat{i}}(w)\beta_{\hat{j}}(w) - \beta_{\hat{i}}(w) c_{\hat{j}}(w))$$ yields 
$$-\frac{1}{z-w} \partial_{j} \partial_{\hat{k}}  P^{\hat{i}\hat{j}}(\gamma) P^{ij}(\gamma) c_i c_{\hat{i}} c_{\hat{j}} b^{\hat{k}} + \frac{1}{z-w} P^{ij}(\gamma) \partial_{i} P^{\hat{i}\hat{j}}(\gamma)  c_{\hat{i}}c_{\hat{j}} \beta_j  - \frac{1}{(z-w)^2} \partial_{j} \partial_{i}  P^{\hat{i}\hat{j}}(\gamma) P^{ij}(\gamma) c_{\hat{i}}c_{\hat{j}} +$$
$$ + \frac{1}{z-w} \partial_{i} \partial_{\hat{k}}  P^{\hat{i}\hat{j}}(\gamma) P^{ij}(\gamma) c_j c_{\hat{i}} c_{\hat{j}} b^{\hat{k}} - \frac{1}{z-w} P^{ij}(\gamma) \partial_{j} P^{\hat{i}\hat{j}}(\gamma)  c_{\hat{i}}c_{\hat{j}} \beta_i  + \frac{1}{(z-w)^2} \partial_{i} \partial_{j}  P^{\hat{i}\hat{j}}(\gamma) P^{ij}(\gamma) c_{\hat{i}}c_{\hat{j}} +$$

$$ + \frac{1}{z-w}  \partial_{\hat{j}}  \partial_{k} P^{ij}(\gamma) P^{\hat{i}\hat{j}}(\gamma) c_i c_{j} b^k c_{\hat{i}} + \frac{1}{z-w} \partial_{\hat{i}} P^{ij}(\gamma)  P^{\hat{i}\hat{j}}(\gamma) c_i c_j \beta_{\hat{j}}  + \frac{1}{(z-w)^2}  \partial_{\hat{j}}  \partial_{\hat{i}} P^{ij}(\gamma) P^{\hat{i}\hat{j}}(\gamma) c_i c_j -$$
$$ - \frac{1}{z-w} \partial_{\hat{i}} \partial_{k} P^{ij}(\gamma) P^{\hat{i}\hat{j}}(\gamma) c_i c_{j} b^k c_{\hat{j}} - \frac{1}{z-w} \partial_{\hat{j}}  P^{ij}(\gamma)  P^{\hat{i}\hat{j}}(\gamma) c_i c_j \beta_{\hat{i}} -  \frac{1}{(z-w)^2} \partial_{\hat{i}}  \partial_{\hat{j}}  P^{ij}(\gamma) P^{\hat{i}\hat{j}}(\gamma) c_i c_j.$$ 

Term $\partial_k P^{ij}(\gamma(z))c_i(z) c_j(z) b^k(z) \partial_{\hat{k}} P^{\hat{i}\hat{j}}(\gamma(w)) c_{\hat{i}}(w)c_{\hat{j}}(w) b^{\hat{k}}(w)$ yields 

$$ \frac{1}{z-w} \partial_{\hat{i}}  P^{ij}(\gamma) \partial_{\hat{k}} P^{\hat{i}\hat{j}}(\gamma)  c_i c_j c_{\hat{j}} b^{\hat{k}} - \frac{1}{z-w}  \partial_{\hat{j}}  P^{ij}(\gamma) \partial_{\hat{k}} P^{\hat{i}\hat{j}}(\gamma) c_i c_j c_{\hat{i}} b^{\hat{k}} +$$ $$+ \frac{1}{(z-w)^2} \partial_{\hat{i}}  P^{ij}(\gamma) \partial_{i} P^{\hat{i}\hat{j}}(\gamma) c_j c_{\hat{j}} -  \frac{1}{(z-w)^2} \partial_{\hat{i}}  P^{ij}(\gamma) \partial_j P^{\hat{i}\hat{j}}(\gamma) c_i c_{\hat{j}}  +  $$ 

$$+ \frac{1}{z-w}  \partial_k P^{ij}(\gamma^k) \partial_{i} P^{\hat{i}\hat{j}}(\gamma)  c_j c_{\hat{i}} c_{\hat{j}} b^{k} 
- \frac{1}{z-w}  \partial_k P^{ij}(\gamma) \partial_{j} P^{\hat{i}\hat{j}}(\gamma)  c_i c_{\hat{i}} c_{\hat{j}} b^{k} -$$
$$- \frac{1}{(z-w)^2} \partial_{\hat{j}}  P^{ij}(\gamma) \partial_{i} P^{\hat{i}\hat{j}}(\gamma) c_j c_{\hat{i}} +  \frac{1}{(z-w)^2} \partial_{\hat{j}}  P^{ij}(\gamma) \partial_{j} P^{\hat{i}\hat{j}}(\gamma) c_i c_{\hat{i}}.$$ 

Grouping up, further we obtain:

Terms of type $c c \beta$:
$$\frac{1}{z-w} P^{\hat{i}\hat{j}}(\gamma) \partial_{\hat{j}}  P^{ij}(\gamma) c_i c_{\hat{i}} \beta_j  - \frac{1}{z-w} P^{\hat{i}\hat{j}}(\gamma)  \partial_{\hat{j}}  P^{ij}(\gamma)  c_i c_j \beta_{\hat{i}}  - \frac{1}{z-w} P^{\hat{i}\hat{j}}(\gamma)  \partial_{\hat{j}}  P^{ij}(\gamma)  c_j c_{\hat{i}} \beta_i  -$$
$$- \frac{1}{z-w} P^{\hat{i}\hat{j}}(\gamma)  \partial_{\hat{i}} P^{ij}(\gamma) c_i c_{\hat{j}} \beta_j + \frac{1}{z-w} P^{\hat{i}\hat{j}}(\gamma) \partial_{\hat{i}} P^{ij}(\gamma)  c_i c_j \beta_{\hat{j}}   + \frac{1}{z-w} P^{\hat{i}\hat{j}}(\gamma) \partial_{\hat{i}}  P^{ij}(\gamma)  c_j c_{\hat{j}} \beta_i+$$
$$- \frac{1}{z-w} P^{ij}(\gamma) \partial_{j} P^{\hat{i}\hat{j}}(\gamma) c_i c_{\hat{i}}\beta_{\hat{j}} + \frac{1}{z-w}P^{ij}(\gamma)  \partial_{j}P^{\hat{i}\hat{j}}(\gamma)  c_i c_{\hat{j}}\beta_{\hat{i}}  - \frac{1}{z-w} P^{ij}(\gamma) \partial_{j} P^{\hat{i}\hat{j}}(\gamma)  c_{\hat{i}}c_{\hat{j}} \beta_i +$$
$$+\frac{1}{z-w} P^{ij}(\gamma)  \partial_{i} P^{\hat{i}\hat{j}}(\gamma) c_j c_{\hat{i}}\beta_{\hat{j}} - \frac{1}{z-w} P^{ij}(\gamma)  \partial_{i} P^{\hat{i}\hat{j}}(\gamma)  c_j c_{\hat{j}}\beta_{\hat{i}} + \frac{1}{z-w} P^{ij}(\gamma)  \partial_{i} P^{\hat{i}\hat{j}}(\gamma)  c_{\hat{i}}c_{\hat{j}} \beta_j,$$ where, as  $P$ is a skew-symmetric, i.e. $P^{ij}(\gamma) = - P^{ji}(\gamma)$, under changing variables $i \leftrightarrow j$ first row is equal to second row and third one -- to forth one.

Rearranging indexes accordingly we can see that second and forth lines are nothing but Jacobi identity:
$$- \frac{1}{z-w} P^{\hat{i}\hat{j}}(\gamma)  \partial_{\hat{i}} P^{ij}(\gamma) c_i c_{\hat{j}} \beta_j - \frac{1}{z-w} P^{\hat{i}j}(\gamma) \partial_{\hat{i}} P^{\hat{j}i}(\gamma)  c_i c_{\hat{j}} \beta_j  - \frac{1}{z-w} P^{\hat{i}i}(\gamma) \partial_{\hat{i}}  P^{j\hat{j}}(\gamma) c_i c_{\hat{j}} \beta_j,$$
$$+\frac{1}{z-w} P^{ij}(\gamma)  \partial_{i} P^{\hat{i}\hat{j}}(\gamma) c_j c_{\hat{i}}\beta_{\hat{j}} + \frac{1}{z-w}  P^{i\hat{i}}(\gamma)  \partial_{i} P^{\hat{j}j}(\gamma)  c_j c_{\hat{i}} \beta_{\hat{j}} + \frac{1}{z-w} P^{i\hat{j}}(\gamma)  \partial_{i} P^{j\hat{i}}(\gamma) c_j c_{\hat{i}}\beta_{\hat{j}}.$$

Terms of type $c c c b$:
$$-\frac{1}{z-w} \partial_{j} \partial_{\hat{k}}  P^{\hat{i}\hat{j}}(\gamma) P^{ij}(\gamma) c_i c_{\hat{i}} c_{\hat{j}} b^{\hat{k}}  + \frac{1}{z-w} \partial_{i} \partial_{\hat{k}}  P^{\hat{i}\hat{j}}(\gamma) P^{ij}(\gamma) c_j c_{\hat{i}} c_{\hat{j}} b^{\hat{k}}  +$$
$$ + \frac{1}{z-w}  \partial_{\hat{j}}  \partial_{k} P^{ij}(\gamma) P^{\hat{i}\hat{j}}(\gamma) c_i c_{j} b^k c_{\hat{i}}  - \frac{1}{z-w} \partial_{\hat{i}} \partial_{k} P^{ij}(\gamma) P^{\hat{i}\hat{j}}(\gamma) c_i c_{j} b^k c_{\hat{j}} +$$
$$ +\frac{1}{z-w} \partial_{\hat{i}}  P^{ij}(\gamma) \partial_{\hat{k}} P^{\hat{i}\hat{j}}(\gamma)  c_i c_j c_{\hat{j}} b^{\hat{k}} - \frac{1}{z-w}  \partial_{\hat{j}}  P^{ij}(\gamma) \partial_{\hat{k}} P^{\hat{i}\hat{j}}(\gamma) c_i c_j c_{\hat{i}} b^{\hat{k}} +$$
$$+ \frac{1}{z-w}  \partial_k P^{ij}(\gamma) \partial_{i} P^{\hat{i}\hat{j}}(\gamma)  c_j c_{\hat{i}} c_{\hat{j}} b^{k} 
- \frac{1}{z-w}  \partial_k P^{ij}(\gamma) \partial_{j} P^{\hat{i}\hat{j}}(\gamma)  c_i c_{\hat{i}} c_{\hat{j}} b^{k},$$
where due to antisymmetric variables first line is equal to second one and third line is equal to forth one under changing variables $m \leftrightarrow \hat{m}$. Moreover, first column is equal to the second under changing variables  $ i \leftrightarrow j$, $\hat{i} \leftrightarrow \hat{j}$. Thus, we get
$$\frac{4}{z-w} \partial_{i} \partial_{\hat{k}}  P^{\hat{i}\hat{j}}(\gamma) P^{ij}(\gamma) c_j c_{\hat{i}} c_{\hat{j}} b^{\hat{k}}  + \frac{4}{z-w}  \partial_k P^{ij}(\gamma) \partial_{i} P^{\hat{i}\hat{j}}(\gamma)  c_j c_{\hat{i}} c_{\hat{j}} b^{k}.$$
But the last is derivative of Jacobi identity. Thus, it is zero.

Terms of type $c c$ do not also survive. Indeed, in OPE $P^{ij}(\gamma) (c_i\beta_j - \beta_i c_j) \partial_{\hat{k}} P^{\hat{i}\hat{j}}(\gamma) c_{\hat{i}}c_{\hat{j}} b^{\hat{k}} + \partial_k P^{ij}(\gamma)c_i c_j b^k P^{\hat{i}\hat{j}}(\gamma) (c_{\hat{i}}\beta_{\hat{j}} - \beta_{\hat{i}} c_{\hat{j}})$ they reduce one another due to symmetry of second derivatives. Moreover, corresponding terms obtained from OPEs $P^{ij}(\gamma) (c_i\beta_j - \beta_i c_j) P^{\hat{i}\hat{j}}(\gamma) (c_{\hat{i}}\beta_{\hat{j}} - \beta_{\hat{i}} c_{\hat{j}})$ and $\partial_k P^{ij}(\gamma)c_i c_j b^k \partial_{\hat{k}} P^{\hat{i}\hat{j}}(\gamma) c_{\hat{i}}c_{\hat{j}} b^{\hat{k}}$ reduce each other as well$._\bigtriangledown$

\begin{Remark} $d_{qL}$ is globally well-defined vertex operator, which raises the fermionic number by $+1,$ as Poisson tensor and  the $bc-\beta\gamma$ system are globally well-defined. According to $[4]$ the cohomology of the $bc-\beta\gamma$ system with differential $d_{qL}$ is again a vertex algebra.\end{Remark}

\section{Cohomology}

While de Rham cohomology groups of manifolds of “finite type” (e.g. compact manifolds) are of finite dimensions, Lichnerowicz–Poisson cohomology groups may have infinite dimension in general. The problem of determining whether the LP-cohomology space is finite dimensional or not is already a difficult open problem for most Poisson structures.  The quantum version of the differential settles down this issue, as $c-\gamma$ part of $bc-\beta\gamma$ system  always emerges as subcomplex. Moreover, general procedure is the following.

Due to OPEs of the $bc-\beta\gamma$ system differential complex $(C =\bigoplus\limits_{k} C^k,\; d_{qL})$ admits filtration by degree of $\hbar$. The filtration is of finite length and the terms of the spectral sequence can be computed inductively.

For a non-negative integer $n$ define space $V_n$ of all elements, containing order of derivatives of $b^*, c_*, \gamma^*, \beta_*$ less or equal $n$. Then on first page of spectral sequence by degree of $\hbar$  the defined filtration by maximal order of derivatives is compatible with the boundary map $d^{\hbar^1}_{qL}$ due to only single contractions. Accordingly, for the $bc-\beta\gamma$ system $\mathbb{D}[b^{*}, c_{*},  \gamma^{*}, \beta_{*}]$ we have such representation $$\mathbb{D}[b^{*}, c_{*},  \gamma^{*}, \beta_{*}] \cong V_0 \oplus\bigoplus\limits_{n=0}^{\infty} V_{n+1}/V_{n}$$

As $(V_{n}, d^{\hbar^1}_{qL})$ is subcomplex of $(V_{n+1}, d^{\hbar^1}_{qL})$ and there exists the short exact sequence
$$0 \longrightarrow V_{n} \longrightarrow V_{n+1} \longrightarrow V_{n+1}/V_{n} \longrightarrow 0,$$
it is sufficient to compute cohomology of $(V_{n}, d^{\hbar^1}_{qL})$ for any non-negative integer $n$.

In addition $$V_{n+1} \cong (V_n \otimes \mathbb{R}[\partial^{n+1} c^{*}, \partial^{n+1} \gamma^{*}]) \otimes \mathbb{R}[\partial^{n+1} b^{*}, \partial^{n+1} \beta^{*}],$$
where $(V_n \otimes \mathbb{R}[\partial^{n+1} c^{*}, \partial^{n+1} \gamma^{*}], d^{\hbar^1}_{qL})$ is subcomplex of $(V_{n+1}, d^{\hbar^1}_{qL})$. Thus, mathematical induction can be used. Notice also that differential $d^{\hbar^1}_{qL}$ satisfies Leibniz rule  and it raises the multiplicity degree by $+2$.

It is worthy to point out here that the complex $(\mathbb{R}[ c_{*}, \gamma^{*}] , d^{\hbar^1}_{qL})$ is the classical LP-complex.

To illustrate the machinery we will consider in detail one of quadratic Poisson structures on $\mathbb{R}^2$. All quadratic Poisson structures on $\mathbb{R}^2$ were classified  (their LP-cohomologies were also determined) in $[9]$:
\begin{align*} 
& P_1=\partial_x \wedge \partial_y,\;\;\; \\ & P_2=xy \; \partial_x \wedge \partial_y,\;\;\; \\ & P_3= (x^2+y^2) \; \partial_x \wedge \partial_y, \;\;\; \\ & P_4=y^2 \; \partial_x \wedge \partial_y.
\end{align*} 

Corresponding quantum Lichnerowicz differentials are
\begin{align*} 
& d_{qL1} = \{\oint (c_1\beta_2-\beta_1c_2)dz, -\}, \\ & d_{qL2} = \{\oint \gamma^1 \gamma^2 (c_1\beta_2-\beta_1c_2) dz + (\gamma^2 b^1 +\gamma^1 b^2) c_1 c_2 dz, -\},\\ & d_{qL3} = \{\oint ((\gamma^1)^2 +(\gamma^2)^2) (c_1\beta_2-\beta_1c_2)dz + 2(\gamma^2 b^1 +\gamma^1 b^2) c_1 c_2 dz, -\}, \\ & d_{qL4} = \{\oint (\gamma^2)^2 (c_1\beta_2-\beta_1c_2)dz + 2\gamma^2 b^1 c_1 c_2 dz, -\}.
\end{align*}

First case is symplectic case: it is usual de Rham cohomology. We have only one non-zero cohomology group, namely $H^0 \cong \mathbb{R}$, constants,  since  $d^{\hbar^1}_{qL1}$ is monomial preserving and the complex can be represented as tensor product of $c-\gamma$ and $b-\beta$ parts.

For second case as well as for third one the LP-cohomology space is finite dimensional. Two cases are similar. We will consider the second case  closely.

The LP-cohomology space is of infinite dimension for the last case. Indeed, $c_2 \partial c_2 \partial^2 c_2 \dots \partial^{m+1} c_2$, $(\gamma^1)^m c_2$ and $(\beta_1)^m$ , where $m \in \mathbb{Z}_{\geq 0}$, represent part of cohomology classes.

Let's inspect complex $(\mathbb{D}[b^{*}, c_{*},  \gamma^{*}, \beta_{*}], d_{qL2})$ on the first page of spectral sequence by degree of $\hbar$. To begin with,  values of $d^{\hbar^1}_{qL2}$ for all single elements are written below:

$$d^{\hbar^1}_{qL2}(\gamma^1)=\gamma^1\gamma^2c_2,\; d^{\hbar^1}_{qL2}(\gamma^2)=\gamma^1\gamma^2c_1;\;d^{\hbar^1}_{qL2}(c_1)=\gamma^2c_1c_2,\; d^{\hbar^1}_{qL2}(c_2)=\gamma^1c_1c_2;$$

$$d^{\hbar^1}_{qL2}(b^1)=\gamma^1\gamma^2\beta_2+(\gamma^2 b^1 +\gamma^1 b^2) c_2,\;d^{\hbar^1}_{qL2}(b^2)=\gamma^1\gamma^2\beta_1+(\gamma^2 b^1 +\gamma^1 b^2) c_1;$$

$$d^{\hbar^1}_{qL2}(\beta_1)=\gamma^2 (c_1\beta_2-\beta_1c_2) + b^2c_1c_2,\;d^{\hbar^1}_{qL2}(\beta_2)=\gamma^1 (c_1\beta_2-\beta_1c_2) +b^1c_1c_2.$$

The complex $(\mathbb{R}[ c_{*}, \gamma^{*}] , d^{\hbar^1}_{qL2})$ is the classical LP-complex and from $[9]$ we know its cohomology classes: $1$, $\gamma^1c_1$, $\gamma^2c_2$, $c_1c_2$, $\gamma^1\gamma^2c_1c_2$. However, notice that in contrast with $[9]$ we could obtain that directly and it would not be difficult.

The routine calculations are left for the reader, but we are about to highlight the key points.

The complex $(V_0\cong\mathbb{R}[ c_{*}, \gamma^{*}] \otimes \mathbb{R}[ b^{*}, \beta_{*}], d^{\hbar^1}_{qL2})$ is more complicated and twisted.  Leibniz rule $d^{\hbar^1}_{qL2}(fg) = d^{\hbar^1}_{qL2}(f)g + fd^{\hbar^1}_{qL2}(g)$ is useful here. Part of cohomology classes are
$$ c_1c_2\beta^k_{1}\beta^m_{2}, \; c_1c_2\beta^k_{1}\beta^m_{2} (\beta_1 b^1 - \beta_2 b^2), \text{ where } k,\;m \in \mathbb{Z}_{\geq 0}.$$

Next step is $ (V_{1} \cong (V_0 \otimes \mathbb{R}[\partial c^{*}, \partial \gamma^{*}]) \otimes \mathbb{R}[\partial b^{*}, \partial\beta^{*}], d^{\hbar^1}_{qL2}),$ where $(V_0  \otimes \mathbb{R}[\partial c^{*}, \partial \gamma^{*}], d^{\hbar^1}_{qL2})$ is subcomplex of $(V_1, d^{\hbar^1}_{qL2})$. Notice that $\frac{\partial}{\partial z}(\gamma^1c_1)$, $\frac{\partial}{\partial z}(\gamma^2c_2)$, $\frac{\partial}{\partial z}(c_1 c_2)$, $\frac{\partial}{\partial z}(\gamma^1\gamma^2c_1c_2)$  represent part of cohomology classes, see the Remark in Chapter 2. 

Moreover,  element $c_1 c_2 \partial c_1 \partial c_2  \dots \partial^{m+1} c_1 \partial^{m+1} c_2, \text{ where } m \in \mathbb{Z}_{\geq 0},$ determines a cohomology class of complex $(\mathbb{D}[b^{*}, c_{*},  \gamma^{*}, \beta_{*}], d^{\hbar^1}_{qL2})$. Thus, there are infinitely many non-zero cohomology groups.

\section{Chiral de Rham Operator}

There is impossible to expect that the chiral de Rham differential (see $[2]$) is going to commute with constructed operator $d_{qL}$, as it doesn't happen on classical level. That's why we consider another differential, which luckily commutes with $d_{qL}$.

Using the same notions $\alpha^{i}=\gamma^i+b^i dz$, $\theta_{i}=c_i+\beta_i dz$ and concept of  variation we have

\[
     \left\{
                \begin{array}{ll}
                  \delta(\alpha^{i})=d(\gamma^i)=\partial_z \gamma^i dz\\
                   \delta(\theta_{i})=d(c_i)=\partial_z c_i dz
                \end{array}
              \right.
  \]
  
In other words, 
$$\delta(b^i)=\partial_z \gamma^i, \delta(\gamma^i) =0,\;$$ $$\delta(\beta_i)=\partial_z c_i, \delta(c_i)=0.$$

Define the associated non-zero OPEs for $bc$ and $\beta\gamma$ systems to be (changing sign of $\beta\gamma$ system from usual one):
$$b(z)c(w) \sim \frac{1}{z-w}, c(z)b(w) \sim \frac{1}{z-w};\;$$ $$\beta(z)\gamma(w) \sim \frac{1}{z-w}, \gamma(z)\beta(w) \sim \frac{-1}{z-w}.$$

Then $\partial_z\gamma(z)\beta(w) \sim \frac{1}{(z-w)^2}$ and the above operator of variation is equal to the following operator
$$\delta_{dR}= \{\oint (\partial_z \gamma^i c_i) dz, -\}.$$

\begin{Proposition} $\delta_{dR}^2=0, \; [\delta_{dR}, d_{qL}]=0.$\end{Proposition}

\textit{Proof.} First statement that $\delta_{dR}$ is differential is obvious as there are no singular terms. 

Second one is not harder. We obtain integral of total derivative, indeed:

Term $\partial_k P^{ij}(\gamma(z))c_i(z) c_j(z) b^k(z) \partial_w \gamma^{\hat{i}}(w) c_{\hat{i}}(w) $ yields $\frac{d}{dz}(P^{ij}(\gamma))c_i c_j$. 

Term $P^{ij}(\gamma(z)) (c_i(z)\beta_j(z) - \beta_i(z) c_j(z)) \partial_w \gamma^{\hat{i}}(w) c_{\hat{i}}(w)$ yields $P^{ij}(\gamma) \frac{d}{dz}(c_ic_j)._\bigtriangledown$

\section{Appendix: Nambu-Poisson bracket}

Let us come back to the concept of Poisson manifold. Let M be a smooth finite dimensional manifold and $C^{\infty}(M)$ be the algebra of smooth functions on this manifold. 

A bilinear mapping $\{\ast, \ast \}: C^{\infty}(M) \times C^{\infty}(M) \rightarrow C^{\infty}(M)$ is said to be a Poisson bracket if for any smooth functions $f, g, h \in C^{\infty}(M)$ it satisfies

\;\;\;\; i) $\{f, g \} = - \{g, f \}$ (skew-symmetry);

\;\;\;\; ii) $\{fg, h \} = f \{g, h \}+g \{f, h \}$ (Leibniz rule);

\;\;\;\; iii) $\{f, \{g,h \} \} + \{g, \{h, f \} \} + \{h, \{f, g \} \} = 0$ (Jacobi identity).

For instance consider the 2-dimensional space $\mathbb{R}^2$ with coordinates denoted by $p, q$ and define the bracket by the formula

\[
\{f, g \}  = \frac{\partial (f, g)}{\partial (p, q)} =
\begin{vmatrix}
\partial_p f & \partial_q f  \\ 
\partial_p g & \partial_q g\\
\end{vmatrix}
\]

A generalization of Poisson bracket was proposed by Y. Nambu in [10], where he introduced a ternary bracket of three smooth functions $f, g, h$ defined on the three dimensional space $\mathbb{R}^3$, whose coordinates are denoted by $x, y, z$. This ternary bracket is defined with the help of the Jacobian of a mapping
$$ (x, y, z) \longrightarrow (f(x, y, z) , g(x, y, z), h(x, y, z))$$
as follows

\[
\{f, g, h \} = \frac{\partial (f, g, h)}{\partial (x, y, z)} =
\begin{vmatrix}
\partial_x f & \partial_y f & \partial_z f  \\ 
\partial_x g & \partial_y g & \partial_z g \\ 
\partial_x h & \partial_y h & \partial_z h  \\ 
\end{vmatrix}
\]

Evidently this ternary bracket is totally skew-symmetric. It can be also verified
 that it satisfies the Leibniz rule
$$\{gh, f_1, f_2\}=g \{h, f_1, f_2\}+h \{g, f_1, f_2\},$$
and the identity
$$\{g, h, \{ f_1, f_2, f_3\} \} = \{ \{g, h, f_1 \}, f_2, f_3 \}  + \{f_1, \{g, h, f_2\}, f_3 \}  + \{f_1, f_2, \{g, h,  f_3\} \}. $$

This identity is called Filippov-Jacobi identity and its n-ary version is the basic component of a concept of n-Lie algebra proposed by V. T. Filippov in [11]. So, Poisson bracket can be generalized to any number of arguments. A smooth manifold endowed with a n-ary Nambu-Poisson bracket is called a Nambu-Poisson manifold of nth order [12].

As was mentioned, Nambu-Poisson bracket is generalization of Poisson bracket, but there is opposite direction: Nambu-Poisson brackets can be defined inductively (see Proposition 3, [13]).

\begin{Proposition} An $n$-bracket, $n>2$, is Nambu-Poisson if and only if fixing an argument we get an $(n-1)$-Nambu-Poisson bracket$._\bigtriangledown$\end{Proposition}

The widest generalization we shall need is the notion of a strong homotopy Lie (or $L_{\infty}$) algebra, which is well-known in algebraic homotopy theory, where it originated. This is obtained by allowing for a countable family of multilinear antisymmetric operations of all arities $n \geq 1$, constrained by a countable series of generalizations of the Jacobi identity known as the $L_{\infty}$ identities. This notion admits specializations indexed by subsets $S \subseteq \mathbb{N}$ of arities and which are defined by requiring vanishing of all products of arities not belonging to $S$. This leads to the notion of $L_{S}$ algebra.  The case $S = \{ n \}$, when only a single product of arity $n$ is non-vanishing, recovers the notion of  n-Lie algebras. We refer the reader to [14][15] for more details.

Due to Leibniz rule Nambu-Poisson bracket acts on each factor as a vector field, whence it must be of the form
$$\{f_1, f_2, \dots, f_n \} = P(df_1, df_2, \dots, df_n),$$
where $P$ is a field of n-vectors on a smooth manifold M [16]. It is called a Nambu-Poisson tensor. Remember that if we use the same definition for n = 2, we get a Poisson tensor.

The Nambu-Poisson tensor fields were characterized as follows by L. Takhtajan [12]  (see additionally also [16])

\begin{Proposition} The $n$-vector field $P$ is a Nambu-Poisson tensor of order $n$ $(n \geq 3)$ iff the natural components of $P$  with respect to any local coordinate system $x^{a}$ of $M$ satisfy the equalities:
$$\sum_{k=1}^n [P^{b_1 b_2 \dots b_{k-1} u b_{k+1} \dots b_n} P^{v a_2 a_3 \dots a_{n-1} b_k} + P^{b_1 b_2 \dots b_{k-1} v b_{k+1} \dots b_n} P^{u a_2 a_3 \dots a_{n-1} b_k}] =0,$$
$$\sum_{u=1}^n [P^{a_1 a_2 a_3 \dots a_{n-1} u} \partial_u P^{b_1 b_2 \dots \dots b_n}  - \sum_{k=1}^n P^{b_1 b_2 \dots b_{k-1} u b_{k+1} \dots b_n} \partial_u P^{a_1 a_2 a_3 \dots a_{n-1} b_k}] =0._\bigtriangledown$$\end{Proposition}

A Nambu-Poisson tensor field $P$ of an even order $n = 2k$ satisfies the condition $[[P, P]] = 0$, where the operation is again the canonical Schouten bracket [16]. This suggests the study of generalized Poisson structures: the Nambu-Poisson cohomology.

It is possible (using Proposition 4) to extend the concept  of quantum Lichnerowicz differential to a Nambu-Poisson tensor $P$ of an even order $n = 2k$. This way, \textbf{we obtain universalization of Kontsevich's theory} (to a smooth manifold one can associate the Lie algebras of multi-vector fields and multi-differential operators, where one can encode classical data (Poisson structures) and quantum data (star products); relating these two led Kontsevich to his famous formality theorem that establishes the deformation quantization of Poisson manifolds) [17][18] in the most direct and natural way. Thus, it opens the road for comprehensive pursuing of rational homotopy theory [19].

\section{Appendix: Gromov-Witten theory}

The object of interest in Gromov-Witten theory is a holomorphic map $\phi: \Sigma \rightarrow X$ from genus $g$ Riemann surface $\Sigma$ to manifold (or orbifold) $X$. The number of such maps is equivalent to the Gromov-Witten invariant, which exhibits invariance under complex deformations on $X$. It has origin from topological string theory, namely in Witten’s work [20] on integrals in two dimensional gravity with enumerative meaning of counting instantons (non-trivial solutions of equations of motion) on $X$ of topological string. 

There is a formal definition of Gromov-Witten invariant in algebraic geometry, wherein it can be expressed through cohomology classes on Calabi-Yau manifold $X$. When the target space is an orbifold, the cohomology that is involved in the Gromov-Witten invariant theory is called Chen-Ruan cohomology. This is the
type of cohomology that is sufficient [21] for orbifolds, rather than the orbifold de
Rham cohomology (in the sense that this enlarges the orbifold de Rham cohomology by keeping track of the automorphisms that the cohomology classes might have).

Quantum Lichnerowicz – Poisson complex is also connected to many geometrical invariants. As we have seen, finding its cohomologies requires enormous calculations. However, first page of spectral sequence by degree of $\hbar$  is more handleable. For example, we can apply the general Künneth theorem to 
$$V_{n} \otimes_{V_0} \mathbb{R}[b^{*}, c^{*}, \gamma^{*}, \beta^{*},\partial^{n+1} b^{*}, \partial^{n+1} c^{*}, \partial^{n+1} \gamma^{*}, \partial^{n+1} \beta^{*}],$$
where $\otimes_{V_0}$ means tensor product over ring $V_0 = \mathbb{R}[b^{*}, c^{*}, \gamma^{*}, \beta^{*}].$ \textbf{Corresponding cohomologies yield new invariant.}

\end{document}